\documentclass{amsart}

\usepackage{graphicx}

\newlength{\defbaselineskip}
\newcommand{\setlinespacing}[1]%
           {\setlength{\baselineskip}{#1 \defbaselineskip}}

 \newtheorem{thm}{Theorem}[section]
 \newtheorem{cor}[thm]{Corollary}
 \newtheorem{lem}[thm]{Lemma}
 \newtheorem{prop}[thm]{Proposition}
 \theoremstyle{definition}
 \newtheorem{defn}[thm]{Definition}
 \theoremstyle{remark}
 \newtheorem{rem}[thm]{Remark}

 \numberwithin{equation}{section}
 
 \DeclareMathOperator{\Tr}{Tr}
 \DeclareMathOperator{\Trace}{Trace}

\begin{document}

\title[periodic 3-manifolds and Modular Categories]
 {periodic 3-manifolds and Modular Categories}

\author{Khaled qazaqzeh }

\address{Department of Mathematics, Louisiana State university,
 Baton rouge, La,  70803 USA}

\email{qazaqzeh@math.lsu.edu}
\date{05/20/2005}


\def\baselinestretch{.60}
\setlinespacing{.60}

\begin{abstract}
A $p$-periodic 3-manifold is a 3-manifold that admits a
$\mathbb{Z}_{p}$-action whose fixed point set is a circle. We give
a congruence that relates the quantum invariant of a $p$-periodic
3-manifold associated to any modular category over an integrally
closed ground ring and the corresponding quantum invariant of its
orbit space.
\end{abstract}

\maketitle

\section*{Introduction}
Let $p$ be an odd prime, and let $G$ be the finite cyclic group
$\mathbb{Z}_{p}$. We assume that all 3-manifolds are compact and
closed. The quantum invariant of a 3-manifold can be defined using
any modular category.  In ~\cite{G99}, Gilmer was interested in
studying the relation between the $SU(2)$- invariants of
$p$-periodic 3-manifolds and their quotient manifolds. He obtained
a congruence relating these invariants. His result was obtained by
using the trace formula of topological quantum field theory (see
proposition (\ref{p:t})) and studying Gaussian sums. Chbili used
the results about the Jones polynomial and the Kauffman
multi-bracket of $p$-periodic links to obtain a similar result for
rational homology 3-spheres for the $SO(3)$-invariants in
\cite{C2}. Also in \cite{C1}, he gave similar results for the
$SU(3)$ and the MOO-invariants. Moreover in ~\cite{CL}, Chen and
Le generalized the above results for rational homology spheres
using any complex simple Lie algebra. We give similar results for
all 3-manifolds using any modular category over an integrally
closed ground ring. Our proof takes place completely in the
context of modular categories. We use the surgery descriptions of
$p$-periodic 3-manifold and its orbit manifold, obtained in
\cite{PS01}, to prove the result.

 In section [1], we give a brief exposition on how to calculate the
quantum invariant for any 3-manifold from its surgery description.
In section [2], we discuss the $\mathbb{Z}_{p}$-actions on
3-manifolds and the relation between the link that describes a
$p$-periodic 3-manifold  and the link that describes its orbit
manifold. Some formulas and results regarding the value of colored
ribbon graphs under the covariant functor $F$ will be given in
section [3]. Finally in section [4], we state and prove the main
result.

\section{quantum invariants of 3-manifolds}
Fix a strict modular category ($\mathcal{V}$,$\{v_{i}\}_{i\in I}$)
with ground ring $K$ and a rank $\mathcal{D}\in K$.
\subsection{Introduction}
 A result due to Lickorish and Wallace asserts that every closed
 oriented 3-manifold can be obtained by surgery on $S^{3}$ along a framed link.

\subsection{The $\tau$-invariant of closed 3-manifolds}
Let $M$ be a closed oriented 3-manifold obtained by surgery on
$S^{3}$ along a framed link $L$.  The $\tau$-invariant of
$(M,\Omega)$ associated to ($\mathcal{V,D}$) where $\Omega$ is a
colored ribbon graph in $M$ is given by
\begin{equation} \label{e:rtc}
    \tau(M,\Omega) =
    \Delta^{\sigma(L)}\mathcal{D}^{-\sigma(L)-m-1}\{L,\Omega \}.
\end{equation}
Here $\sigma(L)$ is the signature of the linking matrix of the
link $L$, and $m$ is the number of components of $L$, and $\Delta
= \{U^{-}\}$ where $U^{-}$ denotes the diagram for the unknot with
a single double point and writhe -1.

 We use the notation
\begin{align*}
\{L,\Omega \} & = \sum_{\lambda\in col(L)} \{L,\Omega \}_{\lambda}\\
&= \sum_{\lambda\in
    col(L)}\prod_{i=1}^{m}\dim(\lambda(L_{i}))F(\Gamma(L,\lambda)\cup
    \Omega),\\
\end{align*}
where col($L$) is the set of all mappings from the set of
components of $L$ to $I$ (the set of simple objects), and
$\Gamma(L,\lambda)$ is the ribbon graph obtained by coloring the
$i$-th component of $L$ by $V_{\lambda(i)}$. Here $F$ is the
covariant functor defined in \cite[chapter.\,I]{T94} which assigns
to a $\mathcal{V}$-colored ribbon graph in $\mathbb{R}$ an element
of the ground ring. The material of this section is due to Turaev
~\cite{T94}.
\subsection{The $I$-invariant}

We take the definition of the quantum invariant to be as follows
\[
I(M,\Omega) = \mathcal{D}\tau(M,\Omega).
\]
Our result is simpler when expressed using this normalization.

\begin{thm} \label{t:t}\textnormal{[Turaev]}
$\tau(M,\Omega)$ (or $I(M,\Omega)$) is a topological invariant of the pair
$(M,\Omega)$.
\end{thm}
\textbf{Example.} We know that $S^{3}$ is obtained by doing
surgery on the empty link, i.e $I(S^{3})=1$. Also, $S^{3}$ is
obtained by doing surgery along the Hopf link $H$ with framing 0
on both components. Hence, we conclude $\{H\} = \mathcal{D}^{2}$.

\begin{cor}
$\{L,\Omega\}$ is invariant under Kirby sliding.
\end{cor}

This corollary is really major part of proof of Theorem
(\ref{t:t}).
\newline
Finally, the $\tau$-invariant can be recovered in terms of the
TQFT-theory $(V,Z)$ which is a functor from the category
$\mathcal{C}$ whose objects are closed surfaces and 3-manifolds as
its morphisms (the surfaces and 3-manifolds have banded links
sitting inside of them) to the category of $K$-modules and
$K$-linear homomorphisms, where $V$ is the functor on the surfaces
and $Z$ is the functor on 3-manifolds. In fact, the assigned value
of a closed 3-manifold under $Z$ is a scalar multiplication
homomorphism from the base ring to itself and that scalar is the
$\tau$-invariant of that manifold. For more details of TQFT see
(\cite{T94,BHMV}).

\section{Periodic links and periodic 3-manifolds}
Let $M$ be a closed oriented 3-manifold that is a result of
surgery on $S^{3}$ along the framed link $L$.
\begin{defn}
  A framed link $L$ in $S^{3}$ is said to be $p$-periodic if there
exists a $\mathbb{Z}_{p}$-action on $S^{3}$, with a fixed point
set equal to a circle, that maps $L$ to itself under this action
and $L$ is assumed to be disjoint from the circle.
\end{defn}
\begin{defn}
$M$ is said to be $p$-periodic if there is an orientation
preserving $\mathbb{Z}_{p}$-action with fixed point set equal to a
circle, and the action is free outside this circle.
\end{defn}
Now we list the following two results from ~\cite{PS01} that will
be used in later sections.
\begin{thm}\label{t:nonfree}
There is a $\mathbb{Z}_{p}$-action on $M$  with a fixed point set
equal to a circle iff $M$ can be obtained be as a result of
surgery on a $p$-periodic link $L$ and $\mathbb{Z}_{p}$ acts
freely on the set of the components of $L$.
\end{thm}

\begin{figure}[h]
  \includegraphics[width=1.90in]{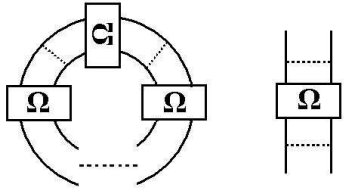}\\
  \caption{Periodic Link and its Quotient}\label{f:figure1}
\end{figure}

By the positive solution of the Smith conjecture, we can represent
any framed $p$-periodic link as a closure of some graph such that
the rotation of this graph about the $z$-axis in $\mathbb{R}^{3}$
(or the circle in $S^{3}$) by 2$\pi/p$ leaves it invariant, i.e $L
= \overline{\Omega}$ (where the bar means the closure of the
graph) see figure ~\ref{f:figure1}.

Let $M_{*} = M/\mathbb{Z}_{p}$ denote the orbit space, then
$M_{*}$ is obtained by surgery on $S^{3}$ along the link $L_{*} =
L/\mathbb{Z}_{p}$ .
\begin{lem}\label{l:basic}
Let $L$ a $p$-periodic link in $S^{3}$. The following are
equivalent
 \begin{enumerate}
  \item $\mathbb{Z}_{p}$ acts freely on the set of components of $L$;
  \item the linking number of each component of the $L_{*}$ the
  axis of the action is congruent to zero modulo $p$;
  \item the number of components of $L$ is  equal to $p$ times the
  number of components of $L_{*}$.
 \end{enumerate}
\end{lem}


\section{Some Results about Traces}
We use two different notions of trace one is the trace of a linear
homomorphism (denoted by $Trace$) in the category of $K$-modules
and the other one is the trace of a ribbon graph (denoted by $Tr$)
in the category of ribbon graphs defined in
~\cite[Chapter.\,1]{T94}
\begin{prop}\label{p:t}
\[
\tau(S^{2}\times S^{1},\overline{\Omega}) =
\Trace_{V(S^{2},l)}(Z_{(S^{2} \times I,\Omega)}),
\]
where $\Omega$ is a colored ribbon $l\times l$ tangle in $S^{2}
\times I$.
\end{prop}
\begin{proof}
This is a special case of the Trace Formula for TQFT
~\cite[Prop.\,1.2]{BHMV} and ~\cite[Ex.\,2.8.1]{T94}.
\end{proof}
\begin{lem}\label{l:newtrace}
\[
\Tr(\Omega)  = \frac{1}{\mathcal{D}^{2}}\sum_{i \in I} \dim(V_{i})
\Trace_{V(S^{2},l+1)}(Z_{(S^{2}\times I,{1_{V_{i}}\otimes
\Omega})}).
\]
\end{lem}
\begin{proof}
Let $H$ stands for the zero-framed Hopf link on both components.
We have
\begin{equation}\label{e:important}
\begin{split}
\Tr(F(\Omega)) & = F(\overline{\Omega}) \qquad \text{by \cite[Cor.\,2.7.2]{T94}}\\
& = \frac{1}{\{H\}}\sum_{\lambda \in col(H)} \{H\}_{\lambda}
F(\overline{\Omega}) \ , \quad \text{as}\ \{H\} = \sum_{\lambda\in col(L)} \{H\}\\
&= \frac{1}{\mathcal{D}^{2}}\sum_{\lambda \in col(H)}
\dim(\lambda)F(\overline{\Omega}\cup H)\ ,\quad\text{as}\ \{H\}
=\mathcal{D}^{2},\
\text{where} \ H \ \text{is unlinked from}\  \overline{\Omega} \\
& =\frac{1}{\mathcal{D}^{2}}\sum_{\lambda \in col(H)}
\dim(\lambda)F(\overline{\Omega^{'}}) \ , \quad
\substack{\text{using the
invariance of}  \ \{L,\Omega\} \  \\
\text{under sliding see figure}\ref{f:figure2}}
\\& = \sum_{i \in I} \dim(V_{i})\mathcal{D}^{-2}\sum_{j \in
I}\dim(V_{j})
F( \overline{\Omega^{'}})\\
& =\sum_{i \in I} \dim(V_{i}) \tau(S^{2}\times
S^{1},\overline{1_{V_{i}}\otimes \Omega})\ ,
\quad \text{by formula (\ref{e:rtc})} \\
& =\sum_{i \in I} \dim(V_{i})
\Trace_{V(S^{2},l+1)}(Z_{(S^{2}\times I,{1_{V_{i}}\otimes
\Omega})})\ , \quad \text{by proposition (\ref{p:t})}
 \end{split}
\end{equation}
\begin{figure}
   \includegraphics[width=3.00in]{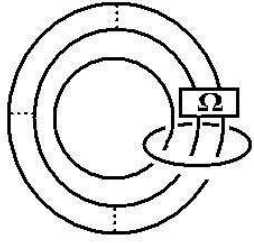}\\
 \caption{The graph $\overline{\Omega^{'}}$ obtained by sliding $\overline{\Omega}$
 over the Hopf link using the second Kirby move }
 \label{f:figure2}
\end{figure}

 \end{proof}
\begin{defn}
Let $J_{p} = (p, \dim(V_{i})^p - \dim(V_{i}))$ be the ideal
generated by $p$ and $\dim(V_{i})^p - \dim(V_{i}), \forall i\in I$
in $K$.
\end{defn}
\begin{cor}\label{c:trace}
Let $\Omega$ be any colored ribbon graph over any modular category
with integrally closed ground ring. Then
\begin{equation}
\Tr(\Omega)^p \equiv \Tr(\Omega^p) \quad \pmod{ J_{p}}.
\end{equation}

\end{cor}
\begin{proof}
If we assume
    \[\dim(V_{i})^{p} = \dim(V_{i}),
\] it follows that
    \[\mathcal{D}^{2p} = \sum_{i \in I} \dim(V_{i})^{2p} = \mathcal{D}^{2}  \quad
    \pmod{p}.
    \]
Hence the result follows from Lemma (\ref{l:newtrace}) and
~\cite[Lem.\,3.5(i)]{CL}, which implies that
\[
\Trace(Z^{p}) \equiv [\Trace(Z)]^{p} \quad \pmod{ p},
\]
where $Z$ is an endomorphism of free $K$-module.
\end{proof}


\section{ Quantum invariants of Periodic 3-manifolds }
 Let $M$ be a 3-manifold that admits a $\mathbb{Z}_{p}$-action
 with a fixed point set equal to a circle. Then we are in
 situation of theorem ~\ref{t:nonfree}. i.e. $M$ is obtained by surgery
 on $S^{3}$ along a framed $p$-periodic link $L$ (see
 figure ~\ref{f:figure1}). We would like to relate the quantum invariant of $M$
 to the quantum invariant of $M_{*} = M/\mathbb{Z}_{p}$. Before
 we do so,  we introduce the following.

\begin{defn}
Let $L$ be a $p$-periodic link, and $\lambda$ be a coloring of
$L$. If $\Gamma(L,\lambda)$ is invariant under the rotation of the
graph that represents $L$ by 2$\pi/p$, then $\lambda$ is called a
$p$-periodic coloring.
\end{defn}
\begin{lem}\label{l:new}
Let $L$ be a $p$-periodic link, such that $L_{*} =
L/\mathbb{Z}_{p}$. Then
\begin{equation}\label{e:e}
\{L\}  \equiv   \{L_{*}\}^{p}   \quad \pmod{J_{p}}.
\end{equation}

\end{lem}
\begin{proof}
Let us start with any coloring of $L$ say $\lambda$, either
$\lambda$ is $p$-periodic or not. Let us assume that $\lambda$ is
not $p$-periodic, i.e $\Gamma(L,\lambda)$ is not invariant under
the rotation by 2$\pi/p$ about the $z$-axis. Hence the $i$-th
rotation of $\Gamma(L,\lambda)$ ( the rotation by $2i\pi/p$)
represents a ribbon graph with the same value under $F$ (since $F$
is an isotopy invariant) and different coloring denoted by
$\lambda_{i}$. So the term with a non-periodic coloring occurs $p$
times. Hence we reduce the summation on the left-hand side to the
periodic colorings. Now the result follows from corollary
(\ref{c:trace}) and the fact that the periodic colorings of $L$
are in one-to-one correspondence with the colorings of $L_{*}$(by
restriction).
\end{proof}
We introduce the notion $\kappa = \Delta\mathcal{D}^{-1}$. Now, we
are ready to give a relation between the quantum invariants of $M$
and $M_{*}$.
\begin{thm}\label{t:main}
Over any modular category with integrally closed ground ring $K$;
we have
\begin{equation}\label{e:main}
I(M) \equiv
 \kappa^{\delta}I(M_{*})^{p} \quad \pmod{ J_{p}},
\end{equation}
for some integer $\delta$.
\end{thm}
\begin{proof} We assume that $M$ and $M_{*}$ are obtained by surgery on $S^{3}$
 along $L$ and $L_{*}$ respectively.

\begin{equation}
\begin{split}
I(M) & = (\Delta\mathcal{D}^{-1})^{\sigma(L)}\mathcal{D}^{-pm}\{L\}\\
& \equiv
(\Delta\mathcal{D}^{-1})^{\sigma(L)}\mathcal{D}^{-pm}\{L_{*}\}^{p}
 \quad \quad \pmod{ J_{p}} \text{by lemma(\ref{l:new})}  \\
& \equiv (\Delta\mathcal{D}^{-1})^{\sigma(L)-p\sigma(L_{*})}
((\Delta\mathcal{D}^{-1})^{\sigma(L_{*})})^{p}
(\mathcal{D}^{-m})^{p} \{L_{*}\}^{p} \quad \pmod{J_{p}} \\
& \equiv \kappa^{\delta}I(M_{*})^{p}  \quad \pmod{ J_{p}}.\\
\end{split}
\end{equation}
 Here $\delta = \sigma(L) - p\sigma(L_{*})$.
\end{proof}
\begin{cor}
\begin{equation}
\tau(M) \equiv \kappa^{\delta}\mathcal{D}^{p-1} \tau(M_{*})^{p}
\quad \pmod{J_{p}},
 \end{equation}
  where $\delta$ and $\kappa$ as defined before.
 \end{cor}
Before we go to the next corollary, we define the total signature
for a knot.
\begin{defn}
Suppose $K$ is a knot in a homology sphere $M_{*}$. Let
$\pi:M\longrightarrow M_{*}$ be the $p$-fold cyclic cover branched
along $K$. It is known that, we can extend this cover to a cover
$W\longrightarrow W_{*}$ of 4-manifolds (where $M=\partial W$ and
$M_{*}=\partial W_{*}$) branching over the surface $Y$. Let $Y
\cdot Y$ denote the self-intersection of $Y$ in $W_{*}$ using the
framing on $K$ obtained from any Seifert surface for $K$ in
$M_{*}$.  In this case, we define the total signature
$\sigma_{p}(K)$
\begin{align*}
\sigma_{p}(K) & = p\sigma(W_{*})-\sigma(W)-\frac{p^{2}-1}{3p}Y\cdot Y \\
& = p\sigma(L_{*})-\sigma(L)-\frac{p^{2}-1}{3p}Y\cdot Y.
\end{align*}
By a well-known argument, using Novikov additivity and the
$G$-signature theory \cite{K}, $\sigma_{p}(K)$ is independent of
the choices made.
\end{defn}
The following corollary generalizes \cite[Th.\,3]{G99}.
\begin{cor} If $M$ is a $p$-fold branched
cyclic cover of a homology sphere along a knot $K$, then
\[
I(M) \equiv \kappa^{-\sigma_{p}(K)} I(M_{*})^{p}   \quad
\pmod{J_{p}},
\] where $\sigma_{p}(K)$ is the total signature of $K$.
\end{cor}
\begin{proof}
The linking matrix of $L$ describes the intersection form on
4-manifold with boundary $M$ which is a branched cover along a
disk with zero self-intersection in a 4-manifold with boundary
$M_{*}$. The corollary now follows by identifying $\delta$ with
the total signature of $K$.
\end{proof}
\begin{cor}
If $M$ is a $p$-fold branched cyclic cover of $S^{3}$ along the
knot $K$, then
\[
I(M) \equiv \kappa^{-\sigma_{p}(K)}   \quad \pmod{J_{p}}.
\]
\end{cor}
\begin{rem}
If $K$ is a knot in $S^{3}$, $\sigma_{p}(K)$ can be identified
with minus the sum of the Tristram-Levine ``p-signatures''.
\end{rem}
\subsection*{Acknowledgment}
The results of this paper were obtained during my Ph.D. studies at
Louisiana State University. I would like to express deep gratitude
to my advisor Dr.Patrick Gilmer whose comments and support were
very important for the completion of this project.

\setlinespacing{1.00}


\begin{thebibliography}{W}
\bibitem[BHMV]{BHMV}
C. ~Blanchet, N. ~Habegger, G. ~Masbaum, and P. ~Vogel,
\emph{Topological Quantum Field Theories Derived from the Kauffman
Bracket}, Topology \textbf{34} (1995), 883-927.
\bibitem[C1]{C1}
N. ~Chbili, \emph{Quantum invariants of finite group actions on
three manifolds},
 Topology and its Applications \textbf{136} (2004), 219-231.
\bibitem[C2]{C2}
N. ~Chbili, \emph{Les invariants $\theta_{p}$ des
3-vari$\acute{e}$t$\acute{e}$s p$\acute{e}$riodiques}, Ann. Inst.
Fourier \textbf{51}(2001) 1135-1150.
\bibitem[CL]{CL}
Q.~Chen, T.~Le, \emph{Quantum invariants of periodic links and
periodic 3-manifolds}, Fund. Math. \textbf{184} (2004), 55--71.
\bibitem[G]{G99}
P.~M. Gilmer, \emph{Quantum invariants of periodic
three-manifolds}, Geometry
  and Topology Monographs \textbf{2} (1999), 157--175.
\bibitem[GB]{GBP02}
P.~M. Gilmer, J.~Kania-Bartoszynska, and J.~H. Przytycki,
\emph{3-manifold
  invariants and periodicity of homology spheres}, Algebraic and Geometric
  Topology \textbf{2} (2002), 825--842.
\bibitem[K]{K}
L. Kauffman, \emph{On Knots}, Number 115, Annals of Mathematics
Studies, Princeton University Press, 1987.
\bibitem[PS]{PS01}
J.~H. Przytycki and M.~V. Sokolov, \emph{Surgeries on periodic
links and
  homology of periodic 3-manifolds}, Math. Proc. Camb. Phil. Soc \textbf{131}
  (2001), 295--307.

\bibitem[T]{T94}
V.~G. Turaev, \emph{Quantum invariants of knots and 3-manifolds},
de Gruyter
  Stud. Math.\textbf{18}, Walter de Gruyter, Berlin; New York, 1994.




\end{thebibliography}
\end{document}